\documentstyle[12pt]{article}
\newcommand{\EBF}[1]{#1&=&}
\newcommand{\LD}{\left.}
\newcommand{\LR}{\left(}
\newcommand{\LC}{\left\{}
\newcommand{\RD}{\right.}
\newcommand{\RRD}{\RD\RD}
\newcommand{\RR}{\right)}
\newcommand{\RC}{\right\}}
\newcommand{\EB}[1]{\\[1mm] #1&=&}
\newcommand{\EC}{\\[0.4mm] & &}
\newcommand{\F}[2]{\frac{#1}{#2}}
\newcommand{\CSTD}[5]{\F{1}{{#1}a^{#2}(\F{n}{2}{#3})_{#4}}\LC(1-a)^{{#5}
            -\F{n}{2}}\LR}
\newcommand{\CSTDA}[6]{\F{1}{{#1}a^{#2}(\F{n}{2}{#3})_{#4}}\LC{#5}(1-a)^{{#6}
            -\F{n}{2}}\LR}
\newcommand{\CSTDB}[4]{\F{a}{{#1}(\F{n}{2}{#2})_{#3}}\LC(1-a)^{{#4}
            -\F{n}{2}}\LR}

\def\be{\begin{equation}}
\def\ee{\end{equation}}
\def\ba{\begin{eqnarray}}
\def\ea{\end{eqnarray}}
\def\ban{\begin{eqnarray*}}
\def\ean{\end{eqnarray*}}
\def\square#1{\mathop{\mkern0.5\thinmuskip\vbox{\hrule
    \hbox{\vrule\hskip#1\vrule height#1 width 0pt\vrule}\hrule}
    \mkern0.5\thinmuskip}}
\def\Square{\mathchoice{\square{8pt}}{\square{7pt}}{\square{6pt}}{\square{5pt}}}

\begin{document}
\thispagestyle{empty}

\title{Complete Computation of DeWitt-Seeley-Gilkey Coefficient $E_4$
 for Nonminimal Operator\\ on Curved Manifolds}

\author{V. P. Gusynin$^1$, V. V. Kornyak$^2$}

\date{}
\maketitle

\centerline{${}^1$Bogolyubov Institute for Theoretical Physics, Kiev, Ukraine}
\centerline{${}^2$Joint Institute for Nuclear Research, Dubna, Russia}

\thispagestyle{empty}
\begin{abstract}
\noindent
Asymptotic heat kernel expansion for nonminimal differential operators on
curved manifolds in the presence of gauge fields is considered. The
complete expressions for the fourth coefficient $(E_4)$ in the heat kernel
expansion for such operators are presented for the first time.
The expressions were computed for general case of manifolds of arbitrary
dimension $n$ and also for the most important case $n=4$. The
calculations have been carried out on PC with the help of a program
written in C.
\end{abstract}

\section{Introduction}
In this paper we continue a computer study of the heat kernel
expansion of elliptic differential operators on curved manifolds and
in the presence of arbitrary gauge fields started in \cite{Ober,JSC}.
The coefficients in the asymptotic expansion of the diagonal heat
kernel elements, called DeWitt-Seeley-Gilkey (DWSG) coefficients
\cite{DeWitt,Seeley,Gilkey}, are of fundamental importance in quantum
field theory, quantum gravity, spectral geometry and topology of manifolds.
Many quantities of interest (such as the effective action,
Green's functions, anomalies in quantum field theory
\cite{DeWitt,Seeley,Gilkey,anomaly}, the indices of elliptic
operators and the invariants of manifolds in spectral geometry
\cite{Ati,invariant})
are expressed in terms of DWSG coefficients~\cite{Barv}.

One more example: the DWSG coefficients, computed for the Sturm--Liouville
(1-dimensional Schr\"odinger) operator, form {\em the Korteweg--de Vries
hierarchy}~\cite{Schim,AvrSchim,Polterovich}.

As is well known, the heat kernel of a positive elliptic differential
operator $A$ of the order $2r$, acting on a bundle of $k$-tensors whose
base is a compact closed $n$-dimensional manifold $M$, has a short-time
expansion~\cite{Seeley,Gilkey,Ati,Gre} in terms of geometrical invariants
\begin{equation}
\langle x\vert e^{-tA}\vert x\rangle \sim\sum_{m\geq 0}
\rm E_m(x\vert A) t^{\textstyle \frac{m-n}{2r}}, \qquad t\to 0_+\ .
\label{expansion}
\end{equation}
The coefficients $E_m(x\vert A)$ depend on the structure of operator $A$.
Being local covariant values of certain dimensionality, these coefficients
can be expressed in terms of coefficient functions of the operator $A$,
curvature (and torsion, in general case) and their covariant derivatives.
Up to now the most complete results were obtained for the second-order
operator
\be
A = - \Box + X,
\label{mini}
\ee
where $\Box=g^{\mu\nu}D_\mu D_\nu$, $D_\mu$ is a covariant derivative
including generally different connections (affine and spinor connections,
gauge fields), $X$ is a matrix in internal space. The computation is based
on the DeWitt ansatz for heat kernel matrix elements
\be
\langle x\vert e^{-tA}\vert x^\prime\rangle= e^{-\frac{\sigma(x,x^\prime)}
{2t}}D^{1/2}(x,x^\prime)\sum_{m=0}^\infty E_m(x,x^\prime\vert A)t^{2m-n
\over 2},
\label{ansatz}
\ee
where $\sigma(x,x^\prime)$ is a geodesic interval which is equal to
a half of square of a distance along a geodesic curve between the points
$x$ and $x^\prime$;
$D(x,x^\prime)=g^{-1/2}\det[-D_\mu D_{\nu^\prime}\sigma](g^\prime)^{-1/2}$
is the van Vleck--Morette determinant, $g=\det g_{\mu\nu}$.
A kernel of the operator $\exp(-tA)$ satisfies a {\em heat equation}.
Substitution of the ansatz (\ref{ansatz}) into the heat equation
leads to recurrent relations for the coefficients $E_m$.
In the case of the operator (\ref {mini}), the explicit expressions for
coincidence limits
$E_m(x\vert A) \equiv E_m(x,x^\prime)\vert_{x=x^\prime}$
have been computed for $m=0,2,4,6,8$ \cite {DeWitt,Gilkey,Sak,Avr}
(odd coefficients are equal to zero for closed manifolds,
they are reserved, however, for manifolds with a boundary).
The most complete expression for $E_8$ has been obtained by
Avramidi~\cite{Avr} (a particular case of a scalar field with
non-minimal conformal connection is considered in~\cite{Xu,Amst}).
In~\cite{torsion} a generalization for the case of Riemann--Cartan spaces
has been considered.
Note that in the flat space case it is possible to go much farther
in computing the coefficients $E_m$ (more exactly ${\rm tr}E_m$).
There are corresponding results up to $E_{10}$ \cite{Ven}
\footnote{For calculating $E_{12}$ (up to total derivative terms)
see Refs.\cite{Schubert,Ven1}. We are grateful to
C.Schubert for bringing our attention to these works.}.
But even in this case, the limit of computation by hands is, apparently,
reached already.
The natural hope of the further progress is connected with development
of analytical computer programs for computing the DWSG
coefficients~\cite{Ober,JSC,Fulling,Bel,Booth}.
For instance, in~\cite{Bel} the DeWitt algorithm for the operator
 $-\Box + X$ in a flat space is implemented on the base of the
computer algebra systems {\rm REDUCE} and {\rm FORM} and in \cite{Booth}
a {\em Mathematica} program for computing heat kernel coefficients is
presented.

However, the DeWitt method does not apply to higher-order operators and
{\em nonminimal\/} operators whose leading term is not a power of the
Laplace operator (for discussion see, for example, \cite{CarinhasFulling,
GusUkr} and a recent paper \cite{BransonAvramidi}). The simplest example
of nonminimal operator is the Navier-Lam\'e operator of classical elasticity,
$\mu\Delta\vec{V} +(\lambda+\mu)\nabla(\nabla\vec{V})$,
involving a coupling among the components of a vector-valued function
(the Lam\'e constants, $\lambda$ and $\mu$, characterize the material).
H. Weyl was apparently the first who investigated the asymptotics of
spectrum of operators of the like type~\cite{Weyl}.
In recent years nonminimal operators of a similar sort have been
encountered by physicists studying the quantization of gauge and
gravitational fields in arbitrary gauges \cite{Barv}.
For example, the quantization of Yang-Mills field in an arbitrary
covariant background gauge leads to the operator
\be
A^{ab}_{\mu\nu}= -\delta_{\mu\nu}\Box^{ab} -
 (\frac{1}{\alpha}-1)D^{ac}_{\mu}
D^{cb}_{\nu} - 2f^{acb}G^c_{\mu\nu} ,
\label{YMoper}
\ee
where $\Box\equiv D_{\mu}D^{\mu}$ is the Laplace operator, $D_{\mu}$
is a covariant derivative containing the external field potential
$A_{\mu}$, $G_{\mu\nu}$ is the field strength and
$f^{abc}$ are the structure constants of a corresponding Lie algebra.
In an analogous way, the quantization of electromagnetic
field in an external gravitational field leads to the operator
\cite{Barth,Cho}.
\be
A_{\mu\nu}=-g_{\mu\nu}\Box - (\frac{1}{\alpha}-1)D_\mu D_\nu + R_{\mu\nu}
\label{EMoper}
\ee
(for an analogous operator in quantum gravity see~\cite{Cho}).
Usually operators (\ref{YMoper}), (\ref{EMoper}) are considered
in the Feynman gauge $\alpha=1$, where they are reduced to
a minimal differential operator (in particular, (\ref{EMoper})
is reduced in this case to the Hodge--de Rham operator
acting on 1-forms), and then are treated by the DeWitt method.
Some works on computation of ${\rm tr}E_m$ are based
on the specific dependencies, e.g. of the operator (\ref{EMoper}), on
the parameter $\alpha$~\cite{Barv,Endo}.
The nontrivial question worth of study is the interplay between
the dependence of a heat kernel on the gauge parameter $\alpha$ and
various invariants appearing in its expansion. This in turn is
important for calculating the gravitational conformal anomaly for
gauge fields in a general covariant gauge and for investigating
possible gauge parameter dependence of the anomaly \cite{Endo}.

The main aim of this paper is to study the heat kernel expansion for
a generic operator with the structure of (\ref{YMoper}),
\be
-g^{\mu\nu}\Box + a D^{\mu}D^{\nu} + X^{\mu\nu},
\label{nmini}
\ee
where $\Box$ is now a covariant generalization of the
Laplace operator, $D^\mu$ is the covariant derivative involving affine
and bundle connections, $X$ is a tensor field and $a$ is a scalar
parameter (bundle indices are assumed implicitly).
The parameter $a$ must satisfy the condition $a<1$ in order to
ensure the positive definiteness of the principal symbol and, hence,
the ellipticity of operator~(\ref{nmini}). The expansion
(\ref{expansion}) is valid for nonminimal operators as well as minimal
ones \cite{Seeley}. In view of inapplicability of the DeWitt method we
follow the approach based on the technique of Widom's covariant
pseudodifferential symbolic calculus~\cite{Widom} developed
in~\cite{Gusynin, GusGor,GGR} (for review see also~\cite{Wincon}).
At present this is the most general method permitting one to handle
operators of general type and of arbitrary order on curved manifolds.

There are several computer implementations of the algorithms for
second-order minimal operators based on general purpose computer algebra
systems (as is mentioned above, in \cite{Bel} REDUCE and FORM and in
\cite{Booth} {\em Mathematica} were used).
However, the calculation of the DWSG coefficients for nonminimal operators
is much more complicated than for second-order minimal ones and is out of
abilities of these programs for nontrivial orders of coefficients.

The algorithm described in~\cite{Gusynin} was implemented as a C
program~\cite{Ober,JSC}. In 1996 we rewrote this program. This last version
allowed us to get for the first time~\cite{GKNIM} the complete
expression of the coefficient $E_4$ for the nonminimal operator (\ref{nmini})
(earlier only partial results for $E_4$ were
known~\cite{Barv,Fradkin,Branson,Pronin}.)

\section{General Outline of Algorithm and Its Implementation}
\vspace*{-0.5pt}
\noindent
In this section we give a short description of the algorithm
and the program for computing the DWSG coefficients.

Let us consider a positive elliptic operator which spectrum lies
inside a contour $C$.
The formula
\be
e^{-tA}=\int_{C}\frac{id\lambda}{2\pi}e^{-t\lambda}(A-\lambda)^{-1},
\ee
allows to express the heat operator $\exp(-tA)$ in terms of the resolvent
$(A-\lambda)^{-1}$ of the operator $A$.
In order to obtain the covariant generalization of the pseudodifferential
calculus method we use the following representation for the matrix elements
of the resolvent
\be
G(x,x',\lambda)\equiv
        \langle x\vert\frac{1}{A-\lambda}\vert x'\rangle
        =\int\frac{d^{n}k}{(2\pi)^{n}\sqrt{g(x')}}e^{il(x,x',k)}
                                        \sigma(x,x',k;\lambda),
\ee
where $\sigma(x,x^\prime,k;\lambda)$ is an amplitude,
$l(x,x^\prime,k)$ is a (real) phase function which is a biscalar
with respect to general coordinate transformations,
$k$ is a wave vector.

The resolvent of an operator $A$ satisfies the equation $(A-\lambda)G=1$
which leads to the equation for the amplitude:
\be
(A(x,D_{a}+iD_{a}l)-\lambda)\sigma(x,x',k;\lambda)=I(x,x'),
\label{ampleq}
\ee
where $I(x,x')$ is a transport function having both bundle and Lorentz
indices.
The covariant generalization of the properties of the flat space phase
and transport functions is reduced to the following relations~\cite{Widom}
\begin{equation}
[\{D_{a_{1}}\ldots D_{a_{m}}\}l]=0,\quad m>1; \qquad
 [\{D_{a_{1}}\ldots D_{a_{m}}\}I]=0,\quad m\geq 1,
\label{lIcond}
\end{equation}
where $\{\ldots\}$ means symmetrizing in all indices,
and $[\ldots]$ means taking the coincidence limit $(x=x')$.
Equations (\ref{lIcond}) together with the ``initial conditions"
$[l]=0, [D_{a}l]=k_{a}$ and $[I]=E$ ($E$ is the unit matrix)
allow one to compute the coincidence limits for  nonsymmetrized
covariant derivatives
$[D_{a_{1}} \ldots D_{a_{m}} l]$ and $[D_{a_{1}} \ldots D_{a_{m}} I]$.
These nonsymmetrized derivatives
are obtained directly from (\ref{lIcond}) by reducing all terms to a
unified index ordering with the help of the Ricci identity.
The resulting expressions are universal polynomials
in the torsion $T^a_{bc}$, curvature tensors $R^a_{bcd}, W_{ab}$ and
their covariant derivatives.
In fact, once computed and stored the coincidence limits
$[D_{a_{1}} \ldots D_{a_{m}} l]$ and $[D_{a_{1}} \ldots D_{a_{m}} I]$
can be used in many calculations for different operators $A$.
The functions $l(x,x^\prime,k)$ and $I(x,x^\prime)$, introduced with
the help of formulas (\ref{lIcond}), play an important role in so-called
{\em intrinsic symbolic calculus} developed in~\cite{Widom}.
In fact, just these universal functions
manifest the geometric properties of a base manifold and a bundle.

Expanding the amplitude $\sigma$ in degrees of homogeneity of
$k: \sigma=\sum_{m=1}^
{\infty}\sigma_{m}(x,x',k;\lambda)$, we obtain the recursion equations for
$\sigma_m$ from Eq.(\ref{ampleq}).
For example, for the operator (\ref{nmini}) these recursion expressions take
the form
\ba
&&A^{ab}\sigma_{0bc}=I^a_c,\nonumber\\
&&A^{ab}\sigma_{1bc}+i\left[-g^{ab}(\Box
l+2D^dlD_d)+a(D^aD^bl+D^alD^b+D^blD^a) \right]\sigma_{0bc}=0,\nonumber\\
&&A^{ab}\sigma_{mbc}+i\left[-g^{ab}(\Box l+2D^dlD_d)+a(D^aD^bl+D^alD^b+
D^blD^a)\right]\sigma_{(m-1)bc}\nonumber\\
&&+(-g^{ab}\Box+aD^aD^b+X^{ab})\sigma_{(m-2)bc}=0,\quad
m\geq 2,
\ea
where the matrix
\ban
A^{ab}=g^{ab}(D^alD_al-\lambda)-aD^alD^bl
\ean
is the principal symbol for the operator (\ref{nmini}). Solving the recursion
equations we obtain $\sigma_m$. The DWSG coefficients are expressed in terms
of integrals of the coincidence limits
$[\sigma_m]$:
\begin{equation}
E_{m}(x\vert A)=\int\frac{d^{n}k}{(2\pi)^{n}\sqrt{g}}
  \int_{C}\frac{id\lambda}{2\pi}e^{-\lambda}[\sigma_{m}](x,k,\lambda)\equiv
  J([\sigma_m]).
\label{eqforem}
\end{equation}
The integrals in (\ref{eqforem}) can be expressed in terms of the Euler
gamma functions and Gauss hypergeometric functions for a wide class of
operators $A$.
The typical integral takes the form
\begin{eqnarray}
\lefteqn{J\Bigl(\frac{k^{2p}k_{a_{1}}\ldots k_{a_{2s}}}
 {(k^{2r}-\lambda)^l[(1-a)k^{2r}-\lambda]^{m}}\Bigr)=} \nonumber\\
 & & g_{\{a_{1}\ldots a_{2s}\}}\frac{\Gamma((p+s+n/2)/r)}
        {(4\pi)^{n/2}2^{s}r\Gamma(n/2+s)\Gamma(l+m)}
        F(m,(p+s+n/2)/r;l+m;a),
\label{integral}
\end{eqnarray}
where $g_{\{a_1\ldots a_{2s}\}}$ is a symmetrized sum of products
of metric tensors.
Using the fact that $m$ and $l$ are integer numbers, one can
express the hypergeometric function in (\ref{integral})
in terms of elementary functions with the help
of the Gauss relation
\begin{equation}
a(1-z)F(a+1,b;c;z)=(c-a)F(a-1,b;c;z)+(2a-c-az+bz)F(a,b;c;z),
\label{Gauss}
\end{equation}
and using then the formula \cite{Prudnik}
\begin{eqnarray}
F(1,b;m;z)=(m-1)!\frac{(-z)^{1-m}}{(1-b)_{m-1}}\Bigl[(1-z)^{m-b-1}-
\sum_{k=0}^{m-2}\frac{(b-m+1)_k}{k!}z^k\Bigr],\nonumber\\
m=1,2,\ldots ,\quad m-b\neq 1,2,\ldots,
\label{Prud}
\end{eqnarray}
where $(a)_k=a(a+1)\ldots(a+k-1)$ is the Pochhammer symbol
(shifted factorial).
During simplification of tensor expressions we use various
symmetry properties of the tensors $R^a_{bcd}, T^a_{bc}, W_{ab}$,
and also the Ricci identity:
\begin{eqnarray}
[D_{a},D_{b}]\varphi^{b_{1}\ldots b_{l}}_{a_{1}\ldots a_{k}}&=&
\sum_{i=1}^{l}R^{b_{i}}_{jab}\varphi^{b_{1}\ldots b_{i-1}j
b_{i+1}\ldots b_{l}}_{a_{1}\ldots a_{k}}
-\sum_{i=1}^{k}R^{j}_{a_{i}ab}\varphi^{b_{1}\ldots b_{l}}_{a_{1}\ldots
a_{i-1}j a_{i+1}\ldots a_{k}}\nonumber
\\
& &+\ T^{j}_{ab}D_{j}\varphi^{b_{1}\ldots b_{l}}_{a_{1}\ldots a_{k}}
+ W_{ab}\varphi^{b_{1}\ldots b_{l}}_{a_{1}\ldots a_{k}},
\end{eqnarray}
the Bianchi identities:
\begin{eqnarray}
D_{a}R^{b}_{cde}+D_{d}R^{b}_{cea}+D_{e}R^{b}_{cad}
  +T^{i}_{ad}R^{b}_{cei}+T^{i}_{de}R^{b}_{cai}+T^{i}_{ea}R^{b}_{cdi}
  =0,\nonumber\\
D_{a}W_{bc}+D_{b}W_{ca}+D_{c}W_{ab}
+W_{ai}T^{i}_{bc}+W_{bi}T^{i}_{ca}+W_{ci}T^{i}_{ab}=0,
\end{eqnarray}
and the cyclic identity:
\begin{equation}
R^{a}_{bcd}+R^{a}_{cdb}+R^{a}_{dbc}
  +D_{b}T^{a}_{cd}+D_{c}T^{a}_{db}+D_{d}T^{a}_{bc}
  +T^{a}_{bi}T^{i}_{cd}+T^{a}_{ci}T^{i}_{db}+T^{a}_{di}T^{i}_{bc}=0
\end{equation}
(in this work we don't consider  torsion, for computing with torsion
see~\cite{Ober,JSC}).

The above algorithm has been implemented in the C language.
The C code of total length about 10000 lines contains about 200
functions for different manipulations with tensors and scalars.
These functions are gathered into two programs COLIM and DWSGCOEF.

The COLIM program computes coincidence limits of the $l$ and $I$
functions and writes them to the disk. Once computed and stored
the coincidence limits, being universal functions,
can be used in many calculations for different operators $A$.

The DWSGCOEF program computes $E_m$ coefficients by the following
steps:\\
1. {\em Reading input information (operator, order $m$, etc.)\/}\\
2. {\em Computing a set of asymptotic operators for constructing
   recursion equations\/}\\
3. {\em Computing $\sigma_m$ with the help of the recursion
   equations\/}\\
4. {\em Taking the coincidence limit $[\sigma_m]$\/}\\
5. {\em Integrating $[\sigma_m]$ to obtain the coefficient $E_m$\/}\\
6. {\em Substituting tensor expressions for
   \([D_{a_{1}} \ldots D_{a_{k}} l]\) and
   \([D_{a_{1}} \ldots D_{a_{k}} I]\) into $E_m$\/}\\
7. {\em Reducing hypergeometric functions to elementary functions in the
   nonminimal case\/}\\
8. {\em Output $E_m$ (and its Lorentz trace in the nonminimal case)\/}

To cut down the intermediate swelling we use {\em term-by-term\/}
strategy, i.e., Steps 4-6 are applied consecutively to single terms
of $\sigma_m$ generated at Step 3.

\section{Computation of $E_4$}
\vspace*{-0.5pt}
\noindent
Note that the coefficient $E_4$ is connected intimately with
the dimension of physical space-time $n=4$ due to the fact that
the Atiyah--Singer index of an elliptic operator on a manifold of
dimension $n$ can be expressed in terms of an integral of $E_n$.
We computed the coefficient $E_4$ for operator (\ref{nmini}) on a
Pentium-75 PC. For completeness we give first the expression for $E_2$.
In the previous version of the program \cite{Ober,JSC} the result was
expressed in terms of hypergeometric functions.
The current version contains a built-in algorithm for reducing
hypergeometric functions to elementary functions and
is based on formulas (\ref{Gauss}), (\ref{Prud}).
\be
E_2=(4\pi)^{-\frac{n}{2}}\LC-C_1X^{ab}
-\frac{C_2}{4}\LR X^{ba}+g^{ab}X_i{}^i\RR
-C_3W^{ab}+C_4R^{ab}+C_5g^{ab}R\RC,
\label{trE_2}
\ee
where
\noindent
\begin{eqnarray*}
C_1&=&\frac{1}{8a(\frac{n}{2}-1)_3}\left\{(1-a)^{-\frac{n}{2}}\left(
4+4n-6a-3an\right)-4-4n+6a-3an-2an^2+an^3\right\},
\\
C_2&=&\frac{1}{2a(\frac{n}{2}-1)_3}\left\{(1-a)^{-\frac{n}{2}}\left(
-4+2a+an\right)+4-2a+an\right\},\\
C_3&=&\frac{1}{4a(\frac{n}{2}-1)_2}\left\{(1-a)^{1-\frac{n}{2}}\left(8-
an\right)-8+8a-3an\right\},\\
C_4&=&\frac{1}{24a(\frac{n}{2}-1)_3}\left\{(1-a)^{-\frac{n}{2}}\left(
12n-12a-8an-an^2+2a^2n+a^2n^2\right)\right.\\
& &\left.-12n+12a+8an-5an^2\right\},\\
C_5&=&\frac{1}{48a(\frac{n}{2}-1)_3}\left\{(1-a)^{-\frac{n}{2}}\left(-
24+12a+8an+an^2-2a^2n-a^2n^2\right)\right.\\
& &\left.+24-12a-an^2+an^3\right\}.
\end{eqnarray*}
The full expression for $E_4$ consists of 73 tensor terms with 43
different scalar coefficients\footnote{There are
28 linear dependencies among these coefficients.}
and is given in Appendix.
Here we reproduce only the trace (with respect to Lorentz indices)
${\rm tr_L}E_4$.
The expression for the trace still contains some novelties in
comparison with the result of Branson et al. \cite{Branson} who
computed ${\rm tr} E_4$ in arbitrary dimension of a space but
without gauge field and neglecting the terms with total derivative.

\begin{eqnarray}
\EBF{{\rm tr_L} E_4}
(4 \pi)^{-\F{n}{2}}\LC\RD -C_1\Box X_i{}^i - C_2 D_i D_j X^{ij}
+ C_3 \LR X_i{}^i X_j{}^j + X_{ij}X^{ij} \RR\nonumber
\EC
+ C_4 X_{ij} X^{ji} + C_5 X_{ij} W^{ij} - C_6 W_{ij} X^{ij}
+ C_7 W_{ij} W^{ij} + C_8 R_{ijkl} R^{ijkl}\nonumber
\EC
- C_9 R_{ij} X^{ij} - C_{10}R_{ij} R^{ij} +C_{11}\Box R
+C_{12} R^2-C_{13} R X_i{}^i \LD\RC,
\label{trE_4}
\end{eqnarray}
\noindent
where
\noindent
\begin{eqnarray*}
\EBF{C_1} 
\CSTD{96}{2}{-2}{4}{1}
96-96a-48an+8a^2n+6a^2n^2+a^2n^3
\RR
\RD
\EC
-96(1-a)^2+32a^2n+2a^2n^2-5a^2n^3+a^2n^4
\LD\RC,
\EB{C_{2}} 
\CSTD{48}{2}{-2}{4}{1}
-48n+48a+24an-4a^2n+a^2n^3
\RR
\RD
\EC
+48n-48a-72an+24an^2+48a^2+4a^2n-24a^2n^2+5a^2n^3
\LD\RC,
\EB{C_{3}} 
\CSTD{16}{}{-1}{3}{}
-4+2a+an
\RR
+4-2a+an
\RC,
\EB{C_{4}} 
\CSTD{16}{}{-1}{3}{}
4+4n-6a-3an
\RD
\RR
\EC
-4-4n+6a-3an-2an^2+an^3
\LD\RC,
\EB{C_{5}} 
\CSTD{192}{2}{-2}{5}{-1}
192n+48n^2+576a+192an-
\RRD
\EC
216an^2-24an^3-1152a^2-944a^2n+136a^2n^2+62a^2n^3+
\EC
2a^2n^4+576a^3+664a^3n+86a^3n^2-31a^3n^3-5a^3n^4-104a^4n
\LD\RR
\EC
-62a^4n^2-a^4n^3+2a^4n^4-192n-48n^2-576a-384an+72an^2
\EC
+576a^2+272a^2n-208a^2n^2+10a^2n^3+4a^2n^4-24a^3n+26a^3n^2
\EC
-9a^3n^3+a^3n^4
\LD\RC,
\EB{C_{6}} 
\CSTD{192}{2}{-2}{5}{-1}
-192n-48n^2+384a+816an+96an^2
\RRD
\EC
+12an^3-768a^2-1120a^2n-148a^2n^2-8a^2n^3-2a^2n^4+384a^3
\EC
+584a^3n+150a^3n^2+7a^3n^3-88a^4n-58a^4n^2-5a^4n^3+a^4n^4
\LD\RR
\EC
+192n+48n^2-384a-624an+48an^2+12an^3+384a^2+304a^2n
\EC
-164a^2n^2+8a^2n^3+2a^2n^4-24a^3n+26a^3n^2-9a^3n^3+a^3n^4
\LD\RC,
\EB{C_{7}} 
\CSTD{48}{}{-1}{2}{1}
-96+22an+an^2+2a^2n-a^2n^2
\RR
\RD\EC
+96-96a+26an-3an^2+an^3
\LD\RC,
\EB{C_{8}} 
\F{(1-a)^{2-\F{n}{2}}-16+n}{180},
\EB{C_{9}} 
\CSTD{24}{}{-1}{3}{}
12n-12a-8an-an^2+2a^2n+a^2n^2
\RR
\RD\EC
-12n+12a+8an-5an^2
\LD\RC,
\EB{C_{10}} 
\CSTD{1440}{}{-1}{3}{}
-360n+360a+296an+60an^2+an^3
\RRD
\EC
-112a^2n-60a^2n^2-2a^2n^3-4a^3n+a^3n^3
\LD\RR
\EC
+360n-360a-296an+116an^2-an^3+an^4
\LD\RC,
\EB{C_{11}} 
\CSTD{480}{2}{-2}{4}{1}
480-240n-240a-120an+16a^2n
\RRD
\EC
+26a^2n^2+11a^2n^3+a^2n^4+4a^3n+4a^3n^2-a^3n^3-a^3n^4
\LD\RR
\EC
-480+240n+720a-360an+120an^2-240a^2+104a^2n-70a^2n^2
\EC
+15a^2n^3-5a^2n^4+a^2n^5
\LD\RC,
\EB{C_{12}} 
\CSTD{576}{}{-1}{3}{}
-144+72a+56an+12an^2
\RRD
\EC
+an^3-16a^2n-12a^2n^2-2a^2n^3-4a^3n+a^3n^3
\LD\RR\EC
+144-72a+16an-16an^2-an^3+an^4
\LD\RC,
\EB{C_{13}} 
\CSTD{48}{}{-1}{3}{}
-24+12a+8an+an^2-2a^2n-a^2n^2
\RR
\RD
\EC
+24-12a-an^2+an^3
\LD\RC.
\end{eqnarray*}           
As one can see, the scalar coefficients in (\ref{trE_2}), (\ref{trE_4})
depend in a rather non-trivial way not only on the gauge
parameter $a$, but also on the space dimension $n$. This is a typical
property of non-minimal and higher order operators \cite{Gusynin}.
In the dimension $n=4$ the coefficients $C_i$ in (\ref{trE_4})
take the form\footnote{Transition to $n = 4$ is not a quite trivial
procedure because the denominators contain the factors $n/2-2$
and this requires to handle uncertainties of the type $0/0$.}
\noindent
\begin{eqnarray*}
C_1&=&\frac{2a-1}{6a^2}\ln(1-a)-\frac{8a^2-21a+6}{36a(1-a)},\quad
C_2=\frac{a+4}{6a^2}\ln(1-a)-\frac{13a^2+6a-24}{36a(1-a)},
\\
C_3&=&\frac{a^2}{48(1-a)^2},\qquad C_4=\frac{13a^2-36a+24}{48(1-a)^2},
\\
C_5&=&-\frac{5a+8}{24a^2}\ln(1-a)-\frac{a^5+11a^4+2a^3+74a^2-180a+96}
{288(1-a)^3},\\
C_6&=&-\frac{7a-8}{24a^2}\ln(1-a)-\frac{a^5-17a^4+154a^3-362a^2+324a-96}
{288a(1-a)^3},
\\
C_7&=&\frac{1}{3},\qquad C_8=-\frac{11}{180},\qquad
C_9=\frac{a(4-3a)}{12(1-a)^2},\qquad C_{10}=\frac{23a^2-46a+8}{360(1-a)^2}
\\
C_{11}&=&-\frac{a^2-5a-2}{12a^2}\ln(1-a)-\frac{133a^2-168a-60}{360a(1-a)},
\\
C_{12}&=&-\frac{a^2+4a-8}{144(1-a)^2},\qquad
C_{13}=\frac{3a^2-6a+4}
{24(1-a)^2}.
\end{eqnarray*}
In the case of symmetric (with respect to Lorentz indices)
matrix $X^{ij}$ the expression for ${\rm tr_L}E_4$ coincides with that
obtained in \cite{Barv,Fradkin,Pronin} (the terms with total derivative
were not computed there). Note also that in this case
the logarithmic dependence on the gauge parameter is contained only
in the coefficients at the terms with total derivative.

As to the full expression for $E_4$ given in Appendix, we should note its
rather cumbersome form. It is quite clear that the next coefficients
will be much larger. Thus, the problem of conveying computational results
of such a kind within scientific society becomes rather crucial.
In this connection, it seems the proposition to create a special
electronic archive for results of computations \cite{Christensen}
makes sense. Such an archive would be useful also for checking by other
authors the obtained results and for comparing efficiency of different
algorithms and programs. For this purpose, however, the problems of
identification of a canonical basis for tensor invariants, elaboration of
efficient algorithms of reducing arbitrary tensor expressions to the
canonical form implied by such a basis and converting from one basis
into another one must be solved.

\section{Conclusion}
\vspace*{-0.5pt}
\noindent
The current version of the program computes $E_4$ for the operator
(\ref{nmini}) on a Pentium-75 PC for 4 h 5 min, whereas computation of
$E_2$ (and $E_4$ for minimal operator) takes trivial time
($< 1$ sec)~\footnote{The previous version\cite{Ober,JSC} computed $E_2$
for nonminimal operator for 3 min 15 sec on AT386.}.
The main reason for such a tremendous swelling are tensor submonomials
of the form
$$ (A_0^{-1})_{ai} (A_0^{-1})^i_j \cdots (A_0^{-1})^k_l (A_0^{-1})^l_b.$$
For nonminimal operator  (\ref{nmini})
$A_0^{-1}$ is the matrix $$(A_0^{-1})_{ab} = \frac{1}{k^2 - \lambda}
\left\{g_{ab} + \frac{a\,k_a\,k_b}{(1-a)\,k^2 - \lambda} \right\},$$
whereas for a minimal operator
$A_0^{-1} = \frac{1}{k^{2r} - \lambda}$,
i.e., a single scalar.
However we can cope in part with this swelling of intermediate expressions
by writing
$(A_0^{-1})_{ab}$ in terms of {\em projectors\/} \cite{Fulling}
$P_{1ab} = g_{ab} - \frac{k_a\,k_b}{k^2}$ and $P_{2ab} =
\frac{k_a\,k_b}{k^2}.$
Using then the standard properties of projectors, $P_i^n=P_i$, $P_iP_j=0,
i\neq j$, we can considerably decrease the computing time.

Besides, some additional work on structuring large output expressions
is needed. Indeed, to be understandable or usable, the result of the
calculation must be put into some kind of standard form.
As one can see, any coefficient  $E_m$ in the expansion (\ref{expansion})
is a linear combination of monomials of a given mass dimensionality.
These monomials are constructed from tensors by multiplication and
index contraction.
Factors in the monomials could be, in  general case, noncommutative
and possess symmetries leading to nonobvious linear dependencies among
the monomials. Thus, we need, first, to identify the linearly independent
basis of tensorial and invariant functionals of curvature, torsion, gauge
fields, and other tensors and their derivatives.

Taking into account additional identities like cyclic, Bianchi etc.
makes the problem of the choice of a canonical basis rather nontrivial.
A certain progress in this direction has been reached
in \cite{FullingKing,Mueller}. For example, in
\cite{FullingKing} a partial basis (containing the Riemann curvature tensor
and its covariant derivatives) for the coefficient $E_8$ for the operator
$A= -\Box + X$ has been constructed.
In  \cite{Mueller} an algorithm for constructing a basis for $E_m$ for
the mentioned operator in a flat space has been proposed.
The analogous problem of construction of a canonical basis for nonminimal
operator (\ref{nmini}) is still waiting for its solution.
Certainly, the choice of a basis can not be unique and different bases
are useful for different purposes.

The next important task is the reduction of an arbitrary tensor expression of
a given mass dimensionality to a chosen canonical basis for $E_m$.
A development of an efficient algorithm for such reduction would be an
essential contribution on the way to full automation of computation of the
DWSG coefficients.
\vspace*{-1mm}

\noindent
We would like to thank S.A. Fulling from Texas A\&M University (College Station,
USA) for fruitful discussions.
This work was supported by Grant No. INTAS-93-2058-ext
``East-West Network in Constrained Dynamical Systems".

\section*{\bf Appendix}
\vspace{-0.5pt}
In this Appendix we give the full expression for the coefficient $E_4$
in expansion (\ref{expansion}) for nonminimal operator (\ref{nmini})
for the cases of manifold of arbitrary dimension $n$ and of dimension $n=4$.
{\footnotesize
\begin{eqnarray*}
\EBF{E_4}
(4\pi)^{-\F{n}{2}}\LC-C_1\Square X^{ab}-C_2\LR
\Square X^{ba}+g^{ab}\Square X_i{}^i+2D_iD^aX^{bi}+2D_iD^bX^{ia}
+2D^aD^bX_i{}^i\RD\RD
\EC
\LD+2g^{ab}D_iD_jX^{ij}\RR
-C_3\LR D_iD^aX^{ib}+D_iD^bX^{ai}+\F{1}{2}X_i{}^bW^{ia}
-\F{1}{2}W_i{}^bX^{ai}\RR+C_4\LR X_i{}^iX^{ab}\RD
\EC
\LD+X_i{}^aX^{ib}+X^a{}_iX^{bi}+
X^{ab}X_i{}^i+X^b{}_iX^{ia}+g^{ab}X_{ij}X^{ji}\RR
+C_5\LR X_i{}^iX^{ba}+X_i{}^aX^{bi}+X_i{}^bX^{ia}\RD
\EC
+X_i{}^bX^{ai}
+X^b{}_iX^{ai}
\LD+X^{ba}X_i{}^i+g^{ab}X_i{}^iX_j{}^j+g^{ab}X_{ij}X^{ij}\RR
+C_6\LR X_i{}^iW^{ab}+X_i{}^aW^{ib}+X^b{}_iW^{ia}\RD
\EC
\LD-W_i{}^aX^{bi}-W_i{}^bX^{ia}\RR
+C_7X^a{}_iX^{ib}+C_8\LR X^a{}_iW^{ib}-W_i{}^aX^{ib}\RR
-C_9\Square W^{ab}-C_{10}W_i{}^aW^{ib}
\EC
-C_{11}W_i{}^bW^{ia}+C_{12}W^{ab}X_i{}^i
\!-C_{13}D_iD_jR^{ijab}+C_{14}D_iD_jR^{iajb}\!-C_{15}R_{ijk}{}^aR^{ijkb}
+C_{16}R_{ijk}{}^aR^{ikjb}
\EC
+C_{17}R_{ij}{}^{ab}X^{ij}
+C_{18}R_{ij}{}^{ab}W^{ij}-C_{19}R_i{}^a{}_j{}^bX^{ij}+C_{20}R_i{}^a{}_j{}^bX^{ji}
+C_{21}R_i{}^a{}_j{}^bW^{ij}-C_{22}\Square R^{ab}
\EC
+C_{23}D_iD^aR^{ib}+C_{24}D_iD^bR^{ia}
-C_{25}R_{ij}R^{ijab}+C_{26}R_{ij}R^{iajb}
-C_{27}\LR R_i{}^aX^{ib}+R_i{}^bX^{ai}\RR
\EC
-C_{28}\LR R_i{}^aX^{bi}+R_i{}^bX^{ia}\RR
-C_{29}\LR R_i{}^aW^{ib}-R_i{}^bW^{ia}\RR+C_{30}R_i{}^aR^{ib}
-C_{31}\LR R^{ab}X_i{}^i+g^{ab}R_{ij}X^{ij}\RR
\EC
+C_{32}D^aD^bR-C_{33}RX^{ab}-C_{34}\LR RX^{ba}+g^{ab}RX_i{}^i\RR
-C_{35}RW^{ab}+C_{36}RR^{ab}+C_{37}\,g^{ab}X_{ij}W^{ij}
\EC
+C_{38}\,g^{ab}W_{ij}X^{ij}
+C_{39}\,g^{ab}W_{ij}W^{ij}+C_{40}\,g^{ab}R_{ijkl}R^{ijkl}
-C_{41}\,g^{ab}R_{ij}R^{ij}+C_{42}\,g^{ab}\Square R
\EC
+C_{43}\,g^{ab}R^2
\LD
\RC.
\end{eqnarray*}
} 
The coefficients $C_i$ in arbitrary dimension $n$:
{\footnotesize
\begin{eqnarray*}
\EBF{C_1} 
\CSTD{192}{3}{-2}{5}{1}
-576-192n+960a+432an+48an^2\!-384a^2-288a^2n-48a^2n^2
\RD\RD
\EC
+8a^3n+6a^3n^2+a^3n^3
\LD\RR
+576+192n-1536a-336an+48an^2+1344a^2+96a^2n-96a^2n^2
\EC
-384a^3+152a^3n+62a^3n^2-27a^3n^3-2a^3n^4+a^3n^5
\LD\RC,
\EB{C_2} 
\CSTD{192}{3}{-2}{5}{1}
192-192a-48an+8a^3n+6a^3n^2+a^3n^3\RR-192+384a-48an\RD
\EC
-192a^2+96a^2n-40a^3n+6a^3n^2+a^3n^3
\LD\RC,
\EB{C_3} 
\CSTD{48}{3}{-2}{5}{1}
-96-48n+96a+24an+48a^2+36a^2n+6a^2n^2-8a^3n-6a^3n^2
\RD\RD
\EC
\LD-a^3n^3\RR
+96+48n-192a-24an+24an^2+48a^2-84a^2n-18a^2n^2+6a^2n^3+48a^3+36a^3n
\EC
-16a^3n^2-3a^3n^3+a^3n^4
\LD\RC,
\EB{C_4} 
\CSTD{64}{2}{-2}{5}{}
-48-24n+80a+36an+4an^2-24a^2-18a^2n-3a^2n^2
\RR\RD
\EC
\LD
+48+24n-80a-12an+8an^2+24a^2-10a^2n-3a^2n^2+a^2n^3
\RC,
\EB{C_5} 
\CSTD{64}{2}{-2}{5}{}
48-48a-12an+8a^2+6a^2n+a^2n^2
\RR
\RD
\EC
\LD
-48+48a-12an-8a^2+6a^2n-a^2n^2
\RC,
\EB{C_6} 
\CSTDA{96}{3}{-2}{5}{2}{1}
-48-96a\!-24an+72a^2\!+54a^2n+9a^2n^2\!-8a^3n-6a^3n^2\!-a^3n^3
\RR\RD
\EC
\LD
+96+96a+96an-336a^2-84a^2n+18a^2n^2+144a^3-4a^3n-12a^3n^2+a^3n^3
\RC,
\EB{C_7} 
\CSTD{64}{2}{-2}{5}{}
-48+72n+24n^2\!-48a-\!156an-36an^2\!+104a^2\!+\!78a^2n+13a^2n^2
\RR\RD
\EC
\LD
+48-\!72n\!-24n^2\!+48a+180an-\!12an^3\!-104a^2\!+22a^2n
+53a^2n^2\!-17a^2n^3\!-3a^2n^4\!+a^2n^5
\RC,
\EB{C_8} 
\CSTD{96}{3}{-2}{5}{1}
96+48n-96a+264an+72an^2-432a^2-372a^2n-90a^2n^2
\RD\RD
\EC\LD
-6a^2n^3+56a^3n+42a^3n^2+7a^3n^3
\RR
-96-48n+192a-264an-96an^2+336a^2+708a^2n
\EC
\LD
+30a^2n^2-36a^2n^3-432a^3-228a^3n+136a^3n^2+15a^3n^3-7a^3n^4
\RC,
\EB{C_9} 
\CSTD{96}{3}{-2}{4}{1}
-768+\!1056a+\!144an\!-\!240a^2\!-\!120a^2n\!-\!8a^3n\!-\!6a^3n^2\!-\!a^3n^3\!-\!4a^4n
\RD\RD
\EC
\LD
+a^4n^3
\RR
\LD
+768\!-\!1824a+240an+1296a^2\!-\!456a^2n+24a^2n^2\!-\!240a^3\!+
208a^3n\!-\!30a^3n^2\!-\!a^3n^3
\RC,
\EB{C_{10}} 
\CSTDA{96}{3}{-2}{4}{3}{1}
128-96a+64an-128a^2-96a^2n-16a^2n^2+36a^3n+20a^3n^2
\RD\RD
\EC\LD
+a^3n^3+4a^4n-a^4n^3
\RR
-384+672a-384an+96a^2+720a^2n-96a^2n^2-384a^3-244a^3n
\EC\LD
+120a^3n^2-11a^3n^3
\RC,
\EB{C_{11}} 
\CSTD{96}{3}{-2}{4}{1}
-384+288a+96an-20a^3n-12a^3n^2-a^3n^3-4a^4n+a^4n^3
\RR\RD
\EC
\LD
+3\LR 128-224a+32an+96a^2-48a^2n+20a^3n-a^3n^3
\RR
\RC,
\EB{C_{12}} 
\CSTD{96}{3}{-2}{5}{1}
96\!-\!384a\!-\!96an+144a^2+108a^2n+18a^2n^2\!-\!8a^3n\!-\!6a^3n^2\!-\!a^3n^3
\RR\RD
\EC
\LD+2\LR
-48+240a+24an-264a^2+6a^2n+9a^2n^2+72a^3-22a^3n-3a^3n^2+a^3n^3\RR
\RC,
\EB{C_{13}} 
\CSTDA{2880}{3}{-2}{4}{2}{1}
-16416+58656a+6072an-29352a^2-16608a^2n-966a^2n^2
\RD\RD
\EC
\LD
+3430a^3n+1863a^3n^2+74a^3n^3+296a^4n-74a^4n^3
\RR
+32832-150144a+4272an
\EC
\LD
+176016a^2-21504a^2n-36a^2n^2-58704a^3
+16108a^3n-738a^3n^2-76a^3n^3+15a^3n^4
\RC,
\EB{C_{14}} 
\CSTD{5760}{2}{-2}{4}{1}
53376+35472n-81840a-54264an-6672an^2+15992a^2n
\RD\RD
\EC
\LD
+8766a^2n^2+385a^2n^3+1540a^3n-385a^3n^3
\RR
-53376-35472n+135216a+63048an
\EC
-11064an^2-81840a^2-15992a^2n+13890a^2n^2-1483a^2n^3
\LD\RC,
\EB{C_{15}} 
\CSTD{552960}{4}{-2}{6}{-2}
4534272+387072n+44390400a+\!20447232an+\!2022912an^2
\RD\RD
\EC
-237126144a^2-98031360a^2n-10807296a^2n^2
-218880a^2n^3+435725568a^3+179327872a^3n
\EC
+21040512a^3n^2+626048a^3n^3+11520a^3n^4
-399736320a^4-167316736a^4n-19284928a^4n^2
\EC
-44256a^4n^3+61168a^4n^4+712a^4n^5
+200829696a^5+90282944a^5n+9158592a^5n^2
\EC
-1199920a^5n^3-237456a^5n^4-7912a^5n^5
-58430208a^6-32038400a^6n-3286544a^6n^2
\EC
+1093416a^6n^3+257132a^6n^4
+14474a^6n^5+10438656a^7+7865248a^7n+1626552a^7n^2
\EC
-189364a^7n^3-103878a^7n^4-9249a^7n^5
-625920a^8-985792a^8n-461792a^8n^2-49696a^8n^3
\EC
\LD
+9512a^8n^4+1508a^8n^5
\RR
+8\LR
-566784-48384n-6682368a-2936064an-277056an^2
\RD
\EC
+16842816a^2+3514080a^2n-564096a^2n^2
-105120a^2n^3-14097696a^3+886960a^3n
\EC
+1234272a^3n^2-13600a^3n^3-20376a^3n^4
+4928832a^4-1609312a^4n-180976a^4n^2
\EC
+126600a^4n^3-1466a^4n^4-2243a^4n^5+90a^4n^6
-1148352a^5+239184a^5n+44488a^5n^2
\EC
-9456a^5n^3-1540a^5n^4+288a^5n^5+45a^5n^6
+78240a^6+91000a^6n+24600a^6n^2-42040a^6n^3
\EC
\LD\LD
-7155a^6n^4+2160a^6n^5+360a^6n^6
\RR\RC,
\EB{C_{16}} 
\CSTD{276480}{4}{-2}{6}{-2}
4534272+387072n+44390400a+\!20447232an+\!2022912an^2
\RD\RD
\EC
-237126144a^2-98031360a^2n-10807296a^2n^2
-218880a^2n^3+435725568a^3+179327872a^3n
\EC
+21040512a^3n^2+626048a^3n^3+11520a^3n^4
-399736320a^4-167298304a^4n-19277248a^4n^2
\EC
-48096a^4n^3+59248a^4n^4+520a^4n^5
+200829696a^5+90209216a^5n+9127872a^5n^2
\EC
-1184560a^5n^3-229776a^5n^4-7144a^5n^5
-58430208a^6-31927808a^6n-3240464a^6n^2
\EC
+1070376a^6n^3+245612a^6n^4+13322a^6n^5
+10438656a^7+7791520a^7n+1595832a^7n^2
\EC
-174004a^7n^3-96198a^7n^4-8481a^7n^5
-625920a^8-967360a^8n-454112a^8n^2-53536a^8n^3
\EC
\LD
+7592a^8n^4+1316a^8n^5
\RR+8
\LR\RD
-566784-48384n-6682368a-2936064an-277056an^2
\EC
+16842816a^2+3514080a^2n-564096a^2n^2
-105120a^2n^3-14097696a^3+886960a^3n
\EC
+1234272a^3n^2-13600a^3n^3-20376a^3n^4
+4928832a^4-1646176a^4n-187696a^4n^2
\EC
+137880a^4n^3+574a^4n^4-2759a^4n^5
-1148352a^5+239184a^5n+44488a^5n^2-9456a^5n^3
\EC
-1540a^5n^4+288a^5n^5+45a^5n^6+78240a^6
+91000a^6n+24600a^6n^2-42040a^6n^3-7155a^6n^4
\EC
\LD\LD
+2160a^6n^5+360a^6n^6
\RR\RC,
\EB{C_{17}} 
\CSTDA{192}{3}{-2}{4}{8}{1}
24+48a-50a^2-16a^2n+6a^2n^2-2a^3-5a^3n-2a^3n^2
\RR\RD
\EC
-192-192a-96an+784a^2-16a^2n-72a^2n^2-384a^3
+224a^3n+80a^3n^2-28a^3n^3
\EC
\LD
-16a^4-36a^4n+10a^4n^2+12a^4n^3-3a^4n^4
\RC,
\EB{C_{18}} 
\CSTD{4608}{3}{-2}{5}{-1}
18432+4608n-55680a-24576an-3456an^2+64512a^2
\RD\RD
\EC
+52992a^2n+10176a^2n^2+96a^2n^3-51840a^3
-48480a^3n-11376a^3n^2-528a^3n^3+24a^3n^4
\EC
+29952a^4+21376a^4n+5680a^4n^2+488a^4n^3
\!-\!16a^4n^4\!-\!5376a^5\!-5680a^5n-\!1756a^5n^2\!-92a^5n^3
\EC
\LD
+19a^5n^4\!+512a^6n+\!320a^6n^2\!+\!16a^6n^3\!-8a^6n^4
\RR
\!+12\LR
-1536-384n+3104a+896an+96an^2
\RD
\EC
-2272a^2-1584a^2n-16a^2n^2+88a^2n^3
+2048a^3+800a^3n-440a^3n^2-84a^3n^3+22a^3n^4
\EC
\LD\LD
+2a^3n^5-448a^4-96a^4n+28a^4n^2
+54a^4n^3-3a^4n^5-192a^5n+60a^5n^3-3a^5n^5
\RR
\RC,
\EB{C_{19}} 
\CSTDA{384}{3}{-2}{5}{4}{1}
-48n\!+\!576a+144an-272a^2\!-\!132a^2n+32a^2n^2\!+12a^2n^3\!-\!16a^3
\RD\RD
\EC
\LD
-68a^3n-44a^3n^2-7a^3n^3
\RR
+192n-2304a-768an+96an^2+3392a^2-48a^2n-464a^2n^2
\EC
-24a^2n^3-1024a^3+864a^3n+408a^3n^2-60a^3n^3
-20a^3n^4-64a^4-160a^4n+4a^4n^2+58a^4n^3
\EC
\LD
-3a^4n^5
\RC,
\EB{C_{20}} 
\CSTDA{384}{3}{-2}{5}{4}{1}
48n+576a+144an-336a^2\!-180a^2n+24a^2n^2+12a^2n^3\!-16a^3
\RD\RD
\EC
\LD
-68a^3n-44a^3n^2-7a^3n^3
\RR
-192n-2304a-384an-96an^2+3648a^2+144a^2n-336a^2n^2
\EC
-72a^2n^3-1280a^3+800a^3n+504a^3n^2-44a^3n^3
-28a^3n^4-64a^4-160a^4n+4a^4n^2+58a^4n^3
\EC
\LD
-3a^4n^5
\RC,
\EB{C_{21}} 
\CSTD{2304}{2}{-2}{5}{-1}
-8832+\!3840n+\!2304n^2\!+\!13824a\!-21888an-\!7296an^2\!-96an^3
\RD\RD
\EC
+14976a^2+30816a^2n+9840a^2n^2+768a^2n^3
-25344a^3-19072a^3n-6544a^3n^2-1064a^3n^3
\EC
-56a^3n^4+5376a^4+5680a^4n+2332a^4n^2+524a^4n^3
+53a^4n^4-128a^5n-224a^5n^2-112a^5n^3
\EC
\LD
-16a^5n^4
\RR
+12\LR
736-320n-192n^2-416a+1872an+256an^2-88an^3-1664a^2-1088a^2n
\RD
\EC
+360a^2n^2+64a^2n^3-20a^2n^4+448a^3+96a^3n
-28a^3n^2-54a^3n^3+3a^3n^5+192a^4n-60a^4n^3
\EC
\LD\LD
+3a^4n^5
\RR
\RC,
\EB{C_{22}} 
\CSTD{11520}{3}{-2}{5}{1}
23040+11520n+167424a+172224an+32592an^2-304320a^2
\RD\RD
\EC
-281616a^2n-78072a^2n^2-6672a^2n^3+62816a^3n
+50288a^3n^2+10234a^3n^3+397a^3n^4
\EC
\LD
+6352a^4n+1588a^4n^2-1588a^4n^3-397a^4n^4
\RR
-23040-11520n-144384a-172224an
\EC
-38352an^2+471744a^2+375888a^2n+24552a^2n^2
-11064a^2n^3-304320a^3-148496a^3n
\EC
\LD
+35536a^3n^2+8270a^3n^3-1375a^3n^4
\RC,
\EB{C_{23}} 
\CSTD{11520}{3}{-2}{5}{1}
187392+41088n-523776a-121152an+2448an^2+125760a^2
\RD\RD
\EC
+87984a^2n+10968a^2n^2\!-792a^2n^3\!+256a^3n+376a^3n^2
\!+170a^3n^3\!+23a^3n^4\!+368a^4n+92a^4n^2
\EC
\LD
-92a^4n^3-23a^4n^4
\RR
+3\LR
-62464-13696n+237056a+22848an-7664an^2-216512a^2
\RD
\EC
\LD\LD
+33200a^2n+\!12968a^2n^2\!-\!1856a^2n^3\!+41920a^3
\!-30160a^3n\!+\!1736a^3n^2\!+\!1802a^3n^3\!-\!263a^3n^4
\RR
\RC,
\EB{C_{24}} 
\CSTDA{5760}{3}{-2}{5}{3}{1}
-31232-8768n+48256a+27488an+3856an^2-29632a^2
\RD\RD
\EC
-27440a^2n-7616a^2n^2-652a^2n^3+5536a^3n
+4496a^3n^2+950a^3n^3+43a^3n^4+688a^4n
\EC
\LD
+172a^4n^2\!-172a^4n^3\!-43a^4n^4
\RR
+93696+26304n-238464a-61920an+1584an^2\!+233664a^2
\EC
\LD
+68976a^2n-1680a^2n^2-540a^2n^3-88896a^3
-26096a^3n+5152a^3n^2+1154a^3n^3-49a^3n^4
\RC,
\EB{C_{25}} 
\CSTD{23040}{4}{-2}{5}{-2}
-61440-15360a-34560an+1645056a^2+379392a^2n
\RD\RD
\EC
-16128a^2n^2-4751616a^3-1306752a^3n
-13632a^3n^2+2400a^3n^3+5413632a^4+2036736a^4n
\EC
+165136a^4n^2+720a^4n^3+476a^4n^4-2663040a^5
-1520160a^5n-258512a^5n^2-17064a^5n^3
\EC
-1384a^5n^4\!+387648a^6\!+\!458544a^6n+\!154720a^6n^2
\!+22692a^6n^3\!+1616a^6n^4\!+46080a^7\!-\!4384a^7n
\EC
\LD
-28600a^7n^2-9692a^7n^3-884a^7n^4-960a^8
-7616a^8n-3284a^8n^2+644a^8n^3+251a^8n^4
\RR
\EC
+4\LR\RD
15360+34560a+16320an-357504a^2-56448a^2n
+10272a^2n^2+438336a^3-4640a^3n
\EC
-18912a^3n^2+2816a^3n^3-119232a^4+38400a^4n
+2524a^4n^2-4152a^4n^3+305a^4n^4+105a^4n^5
\EC
\LD\LD
-11040a^5+3320a^5n+60a^5n^2-170a^5n^3
+30a^5n^4+240a^6-220a^6n+60a^6n^2-5a^6n^3
\RR\RC,
\EB{C_{26}} 
\CSTD{1440}{2}{-2}{4}{1}
-27984-2472n+16248a+9108an+492an^2-1778a^2n
\RD\RD
\EC
\LD
-945a^2n^2-28a^2n^3-112a^3n+28a^3n^3
\RR
+27984+2472n-44232a+2412an+744an^2
\EC
\LD
+16248a^2-4234a^2n-237a^2n^2+91a^2n^3
\RC,
\EB{C_{27}} 
\CSTD{192}{3}{-2}{5}{}
192+96n-384a+72an^2-128a^2-512a^2n-136a^2n^2-4a^2n^3
\RD\RD
\EC
\LD
+416a^3+400a^3n+118a^3n^2+11a^3n^3
-56a^4n-42a^4n^2-7a^4n^3
\RR
-192-96n+384a-96an
\EC
\LD
-120an^2+128a^2+656a^2n+88a^2n^2\!-\!44a^2n^3
\!-\!416a^3\!-\!272a^3n+146a^3n^2+23a^3n^3\!-\!9a^3n^4
\RC,
\EB{C_{28}} 
\CSTD{192}{3}{-2}{5}{}
-192+384a-96an+64a^2+144a^2n+32a^2n^2-160a^3
\RD\RD
\EC
\LD
-144a^3n-38a^3n^2-3a^3n^3+24a^4n+18a^4n^2+3a^4n^3
\RR
+192-384a+192an-64a^2
\EC
\LD
-288a^2n+40a^2n^2+160a^3+48a^3n-34a^3n^2+3a^3n^3
\RC,
\EB{C_{29}} 
\CSTD{96}{3}{-2}{4}{1}
192-288a+96an-144a^2-112a^2n-20a^2n^2+48a^3n+26a^3n^2
\RD\RD
\EC
\LD
+a^3n^3+4a^4n-a^4n^3
\RR
-192+480a-192an-144a^2+400a^2n-52a^2n^2-144a^3-144a^3n
\EC
\LD
+70a^3n^2-7a^3n^3
\RC,
\EB{C_{30}} 
\CSTD{960}{2}{-2}{5}{}
37312+12624n+1184n^2\!-58976a-31944an-5116an^2\!-204an^3
\RD\RD
\EC
+22624a^2+21280a^2n+6142a^2n^2+599a^2n^3
+10a^2n^4-2520a^3n-2050a^3n^2-435a^3n^3
\EC
\LD
-20a^3n^4-160a^4n-40a^4n^2+40a^4n^3+10a^4n^4
\RR
-37312-12624n-1184n^2+58976a
\EC
\LD
+13288an-1196an^2-388an^3-22624a^2
-1120a^2n+2010a^2n^2+85a^2n^3-56a^2n^4
\RC,
\EB{C_{31}} 
\CSTD{96}{2}{-2}{5}{}
-72n+128a+96an+16an^2-80a^2-68a^2n-16a^2n^2-a^2n^3
\RD\RD
\EC
\LD
+8a^3n+6a^3n^2+a^3n^3
\RR
\LD
+2\LR
36n-64a-48an+10an^2+40a^2+2a^2n-7a^2n^2+a^2n^3
\RR
\RC,
\EB{C_{32}} 
\CSTD{5760}{3}{-2}{5}{1}
11520-2880n+100416a+34992an+2472an^2-63072a^2
\RD\RD
\EC
-51240a^2n-10836a^2n^2-492a^2n^3+7432a^3n
+5878a^3n^2+1157a^3n^3+38a^3n^4+608a^4n
\EC
\LD
+152a^4n^2-152a^4n^3-38a^4n^4
\RR
-11520+2880n-88896a-43632an-1032an^2+163488a^2
\EC
\LD
+38904a^2n-6348a^2n^2-384a^2n^3-63072a^3
-1072a^3n+4862a^3n^2-227a^3n^3-41a^3n^4
\RC,
\EB{C_{33}} 
\CSTD{192}{2}{-2}{5}{}
-288-144n+448a+224an+44an^2+4an^3-112a^2-140a^2n
\RD\RD
\EC
\LD
-56a^2n^2-7a^2n^3+24a^3n+18a^3n^2+3a^3n^3
\RR
+288+144n-448a-80an+28an^2-4an^3
\EC
+112a^2+52a^2n+16a^2n^2-17a^2n^3-2a^2n^4+a^2n^5
\LD\RC,
\EB{C_{34}} 
\CSTD{192}{2}{-2}{5}{}
288-320a-96an-4an^2+80a^2+68a^2n+16a^2n^2+a^2n^3
\RD\RD
\EC
\LD\LD
-8a^3n-6a^3n^2-a^3n^3
\RR
-288+320a-48an+4an^2-80a^2+20a^2n-4a^2n^2+a^2n^3
\RC,
\EB{C_{35}} 
\CSTD{96}{2}{-2}{4}{1}
-288+144a+88an+8an^2-24a^2n
-14a^2n^2-a^2n^3
\RD\RD
\EC
\LD\LD
-4a^3n+a^3n^3
\RR
+288-432a+56an-8an^2+144a^2-32a^2n+14a^2n^2-3a^2n^3
\RC,
\EB{C_{36}} 
\CSTD{576}{2}{-2}{5}{}
-432n+576a+528an+144an^2+12an^3-288a^2-392a^2n
\RD\RD
\EC
-176a^2n^2-28a^2n^3-a^2n^4+64a^3n
+64a^3n^2+20a^3n^3+2a^3n^4+16a^4n+4a^4n^2-4a^4n^3
\EC
\LD
-a^4n^4
\RR
\LD
+432n\!-\!576a\!-\!528an+72an^2\!-\!12an^3\!+288a^2
+104a^2n+20a^2n^2+10a^2n^3\!-\!5a^2n^4
\RC,
\EB{C_{37}} 
\CSTD{192}{2}{-2}{5}{-1}
-192+1152a+192an-1680a^2-648a^2n-36a^2n^2+768a^3
\RD\RD
\EC
\LD
+628a^3n+96a^3n^2\!+2a^3n^3\!-24a^4\!-178a^4n
-63a^4n^2\!-5a^4n^3\!-24a^5-2a^5n+9a^5n^2\!+2a^5n^3
\RR
\EC
+192-960a-96an+720a^2+24a^2n-36a^2n^2-48a^3
-20a^3n+24a^3n^2-4a^3n^3-24a^4
\EC
\LD
+26a^4n-9a^4n^2+a^4n^3
\RC,
\EB{C_{38}} 
\CSTDB{192}{-2}{5}{-1}
-624-24n+1056a+228an+12an^2-456a^2-238a^2n-39a^2n^2
\RD\RD
\EC
\LD
-2a^2n^3+24a^3+42a^3n+21a^3n^2+3a^3n^3
\RR
+624+24n-432a+108an+24a^2-26a^2n
\EC
\LD
+9a^2n^2-a^2n^3
\RC,
\EB{C_{39}} 
\CSTD{192}{2}{-2}{4}{1}
576-384a-192an+44a^2n+24a^2n^2+a^2n^3+4a^3n-a^3n^3
\RR\RD
\EC
\LD
-576+960a-96an-384a^2+116a^2n-4a^2n^2-5a^2n^3+a^2n^4
\RC,
\EB{C_{40}} 
\F{(1-a)^{2-\F{n}{2}}+n-5}{180(n-4)},
\EB{C_{41}} 
\CSTD{5760}{2}{-2}{5}{}
2160n-4800a-3600an-600an^2\!+3360a^2\!+2984a^2n+776a^2n^2
\RD\RD
\EC
\LD
+64a^2n^3+a^2n^4-448a^3n-352a^3n^2
-68a^3n^3-2a^3n^4-16a^4n-4a^4n^2+4a^4n^3+a^4n^4
\RR\EC\LD
-2160n+4800a+3600an-480an^2\!-3360a^2-520a^2n
+484a^2n^2\!-54a^2n^3-a^2n^4+a^2n^5
\RC,
\EB{C_{42}} 
\CSTD{960}{3}{-2}{5}{1}
1920-1920a-480an+64a^3n+56a^3n^2\!+14a^3n^3\!+a^3n^4
\RD\RD
\EC\LD
+16a^4n+4a^4n^2-4a^4n^3-a^4n^4
\RR
-1920+3840a-480an-1920a^2+960a^2n-320a^3n
\EC\LD
+64a^3n^2-14a^3n^3-a^3n^4+a^3n^5
\RC,
\EB{C_{43}} 
\CSTD{2304}{2}{-2}{5}{}
1728-2112a-720an-48an^2+672a^2+584a^2n+152a^2n^2
\RD\RD
\EC\LD
+16a^2n^3+a^2n^4-64a^3n-64a^3n^2-20a^3n^3
-2a^3n^4-16a^4n-4a^4n^2+4a^4n^3+a^4n^4
\RR\EC\LD
-1728+2112a-144an+48an^2-672a^2+104a^2n-8a^2n^2-12a^2n^3-a^2n^4+a^2n^5
\RC.
\end{eqnarray*}               } 

The coefficients $C_i$ in the dimension $n=4$:
{\footnotesize
\begin{eqnarray*}
\EBF{C_{1}} 
\F{7-11a+a^2}{24a^3}\ln(1-a)
+\F{84-174a+112a^2-15a^3}{288a^2(1-a)},
\EB{C_{2}} 
\F{-1+a+a^2}{24a^3}\ln(1-a)
+\F{-12+18a+8a^2-7a^3}{288a^2(1-a)},
\EB{C_{3}} 
\F{3+a-2a^2}{12a^3}\ln(1-a)
+\F{36-6a-36a^2-5a^3}{144a^2(1-a)},
\EB{C_4} 
\F{3}{32a^2}\ln(1-a)
+\F{36-54a+12a^2+7a^3}{384a(1-a)^2},
\EB{C_5} 
-\F{1}{32a^2}\ln(1-a)
+\F{-12+18a-4a^2-a^3}{384a(1-a)^2},
\EB{C_6} 
\F{1+5a-4a^2}{24a^3}\ln(1-a)
+\F{12+54a-80a^2+13a^3}{288a^2(1-a)},
\EB{C_7} 
-\F{13}{32a^2}\ln(1-a)
+\F{-156+426a-340a^2+71a^3}{384a(1-a)^2},
\EB{C_8} 
\F{-3-25a+14a^2}{24a^3}\ln(1-a)
+\F{-36-282a+324a^2+53a^3}{288a^2(1-a)},
\EB{C_9} 
\F{16-18a-3a^2+a^3}{12a^3}\ln(1-a)
+\F{96-156a+20a^2+31a^3}{72a^2(1-a)},
\EB{C_{10}} 
\F{-8-18a+30a^2-3a^3}{12a^3}\ln(1-a)
+\F{-24-42a+121a^2-40a^3}{36a^2(1-a)},
\EB{C_{11}} 
\F{8-6a-6a^2+a^3}{12a^3}\ln(1-a)
+\F{24-30a-13a^2+16a^3}{36a^2(1-a)},
\EB{C_{12}} 
\F{-1+7a-2a^2}{24a^3}\ln(1-a)
+\F{-12+90a-64a^2-a^3}{288a^2(1-a)},
\EB{C_{13}} 
\F{684-2772a+1863a^2-148a^3}{360a^3}\ln(1-a)
+\F{1368-6228a+6330a^2-1409a^3}{720a^2(1-a)},
\EB{C_{14}} 
\F{-4068+4383a-385a^2}{720a^2}\ln(1-a)
+\F{-8136+12834a-3407a^2}{1440a(1-a)},
\EB{C_{15}} 
\F{-1584-47624a+31554a^2+6483a^3+2906a^4}{17280a^4}\ln(1-a)
\EC
+(-95040-2524800a+11482440a^2-18335200a^3+12961138a^4
\EC
\F{-3498802a^5-762607a^6+1727374a^7-1316820a^8+353710a^9)}
{1036800a^3(1-a)^4},
\EB{C_{16}} 
\F{-1584-47624a+31554a^2+6483a^3+3098a^4}{8640a^4}\ln(1-a)
\EC
+(-95040-2524800a+11482440a^2-18421600a^3+13306738a^4
\EC
\F{-4017202a^5-417007a^6+1640974a^7-1316820a^8+353710a^9)}
{518400a^3(1-a)^4},
\EB{C_{17}} 
\F{-4-12a-9a^2}{24a^3}\ln(1-a)
+\F{-24-60a-14a^2+83a^3+37a^4}{144a^2(1-a)},
\EB{C_{18}} 
\F{-96+257a-101a^2+16a^3}{288a^3}\ln(1-a)
\EC
+\F{-1152+5964a-10458a^2+5708a^3
+1807a^4-896a^5-1752a^6+864a^7}{3456a^2(1-a)^3},
\EB{C_{19}} 
\F{2-10a-15a^2}{24a^3}\ln(1-a)
+\F{12-66a-62a^2+84a^3+37a^4}{144a^2(1-a)},
\EB{C_{20}} 
\F{-2-14a-15a^2}{24a^3}\ln(1-a)
+\F{-12-78a-46a^2+102a^3+37a^4}{144a^2(1-a)},
\EB{C_{21}} 
\F{-113+173a-40a^2}{144a^2}\ln(1-a)
+\F{-1356+5466a-6428a^2+1937a^3\!-592a^4+1752a^5\!-864a^6}
{1728a(1-a)^3},
\EB{C_{22}} 
\F{-180-3768a+4323a^2-397a^3}{720a^3}\ln(1-a)
+\F{-120-2452a+4158a^2-1179a^3}{480a^2(1-a)},
\EB{C_{23}} 
\F{-916+1608a+39a^2-23a^3}{720a^3}\ln(1-a)
+\F{-5496+12396a-3674a^2-559a^3}{4320a^2(1-a)},
\EB{C_{24}} 
\F{518-1200a+1167a^2-129a^3}{360a^3}\ln(1-a)
+\F{1554-4377a+5042a^2-1712a^3}{1080a^2(1-a)},
\EB{C_{25}} 
\F{80+520a-2182a^2+1549a^3-28a^4}{720a^4}\ln(1-a)
\EC
+\F{1920\!+\!5760a\!-\!87728a^2\!+\!280624a^3\!-\!421348a^4
\!+\!331500a^5\!-\!129802a^6\!+\!18484a^7\!+\!810a^8\!+\!5a^9}
{17280a^3(1-a)^4},
\EB{C_{26}} 
\F{1578-945a+56a^2}{360a^2}\ln(1-a)
+\F{3156-3468a+499a^2}{720a(1-a)},
\EB{C_{27}} 
\F{-3-10a+7a^2}{24a^3}\ln(1-a)
+\F{-36-66a+252a^2-109a^3-38a^4}{288a^2(1-a)^2},
\EB{C_{28}} 
\F{1+2a-3a^2}{24a^3}\ln(1-a)
+\F{12+6a-68a^2+63a^3-10a^4}{288a^2(1-a)^2},
\EB{C_{29}} 
\F{-4-6a+13a^2-a^3}{12a^3}\ln(1-a)
+\F{-24-24a+100a^2-31a^3}{72a^2(1-a)},
\EB{C_{30}} 
\F{-556+355a-20a^2}{120a^2}\ln(1-a)
+\F{-3336+7134a-4327a^2+544a^3}{720a(1-a)^2},
\EB{C_{31}} 
\F{3-2a}{24a^2}\ln(1-a)
+\F{12-26a+16a^2-a^3}{96a(1-a)^2},
\EB{C_{32}} 
\F{-1458+1005a-76a^2}{720a^2}\ln(1-a)
+\F{-2916+3468a-659a^2}{1440a(1-a)},
\EB{C_{33}} 
\F{3-2a}{16a^2}\ln(1-a)
+\F{36-46a+11a^3}{192a(1-a)^2},
\EB{C_{34}} 
\F{-3+2a}{48a^2}\ln(1-a)
+\F{-12+26a-16a^2+3a^3}{192a(1-a)^2},
\EB{C_{35}} 
\F{6-7a+a^2}{12a^2}\ln(1-a)
+\F{12-20a+9a^2}{24a(1-a)},
\EB{C_{36}} 
\F{9-12a+2a^2}{72a^2}\ln(1-a)
+\F{12-34a+36a^2-13a^3}{96a(1-a)^2},
\EB{C_{37}} 
\F{4-28a+5a^2}{96a^2}\ln(1-a)
+\F{48-456a+988a^2-778a^3+182a^4+13a^5-a^6}{1152a(1-a)^3},
\EB{C_{38}} 
\F{5a}{32}\ln(1-a)
+\F{a^2(180-450a+330a^2-57a^3+a^4)}{1152(1-a)^3},
\EB{C_{39}} 
\F{-12+12a-a^2}{24a^2}\ln(1-a)
+\F{-6+10a-3a^2}{12a(1-a)},
\EB{C_{40}} 
\F{2-\ln(1-a)}{360},
\EB{C_{41}} 
\F{-45+60a-2a^2}{720a^2}\ln(1-a)
+\F{-180+526a-452a^2+91a^3}{2880a(1-a)^2},
\EB{C_{42}} 
\F{-5+5a+5a^2-a^3}{60a^3}\ln(1-a)
+\F{-60+90a+64a^2-59a^3}{720a^2(1-a)},
\EB{C_{43}} 
\F{-9+12a-2a^2}{288a^2}\ln(1-a)
+\F{-36+118a-116a^2+37a^3}{1152a(1-a)^2}.
\end{eqnarray*}               } 

\end{document}